\documentclass[twoside,11pt]{article}
\usepackage{amsmath,amsthm,amsfonts,amsopn,amssymb}

\makeatletter

\def\lhead{\@ifnextchar[{\@xlhead}{\@ylhead}}
\def\@xlhead[#1]#2{\gdef\@elhead{#1}\gdef\@olhead{#2}}
\def\@ylhead#1{\gdef\@elhead{#1}\gdef\@olhead{#1}}

\def\chead{\@ifnextchar[{\@xchead}{\@ychead}}
\def\@xchead[#1]#2{\gdef\@echead{#1}\gdef\@ochead{#2}}
\def\@ychead#1{\gdef\@echead{#1}\gdef\@ochead{#1}}

\def\rhead{\@ifnextchar[{\@xrhead}{\@yrhead}}
\def\@xrhead[#1]#2{\gdef\@erhead{#1}\gdef\@orhead{#2}}
\def\@yrhead#1{\gdef\@erhead{#1}\gdef\@orhead{#1}}

\def\lfoot{\@ifnextchar[{\@xlfoot}{\@ylfoot}}
\def\@xlfoot[#1]#2{\gdef\@elfoot{#1}\gdef\@olfoot{#2}}
\def\@ylfoot#1{\gdef\@elfoot{#1}\gdef\@olfoot{#1}}

\def\cfoot{\@ifnextchar[{\@xcfoot}{\@ycfoot}}
\def\@xcfoot[#1]#2{\gdef\@ecfoot{#1}\gdef\@ocfoot{#2}}
\def\@ycfoot#1{\gdef\@ecfoot{#1}\gdef\@ocfoot{#1}}

\def\rfoot{\@ifnextchar[{\@xrfoot}{\@yrfoot}}
\def\@xrfoot[#1]#2{\gdef\@erfoot{#1}\gdef\@orfoot{#2}}
\def\@yrfoot#1{\gdef\@erfoot{#1}\gdef\@orfoot{#1}}

\newdimen\headrulewidth
\newdimen\footrulewidth
\newdimen\plainheadrulewidth
\newdimen\plainfootrulewidth
\newdimen\headwidth
\newif\if@fancyplain \@fancyplainfalse
\def\fancyplain#1#2{\if@fancyplain#1\else#2\fi}

\footrulewidth\z@
\plainheadrulewidth\z@
\plainfootrulewidth\z@

\def\@fancyhead#1#2#3#4#5{#1\hbox to\headwidth{\vbox{\hbox
{\rlap{\parbox[b]{\headwidth}{\raggedright#2\strut}}\hfill
\parbox[b]{\headwidth}{\centering#3\strut}\hfill
\llap{\parbox[b]{\headwidth}{\raggedleft#4\strut}}}\headrule}}#5}

\def\@fancyfoot#1#2#3#4#5{#1\hbox to\headwidth{\vbox{\footrule
\hbox{\rlap{\parbox[t]{\headwidth}{\raggedright#2\strut}}\hfill
\parbox[t]{\headwidth}{\centering#3\strut}\hfill
\llap{\parbox[t]{\headwidth}{\raggedleft#4\strut}}}}}#5}

\def\headrule{{\if@fancyplain\headrulewidth\plainheadrulewidth\fi
\hrule\@height\headrulewidth\@width\headwidth \vskip-\headrulewidth}}

\def\footrule{{\if@fancyplain\footrulewidth\plainfootrulewidth\fi
\vskip-0.3\normalbaselineskip\vskip-\footrulewidth
\hrule\@width\headwidth\@height\footrulewidth\vskip0.3\normalbaselineskip}}

\def\ps@fancy{
\let\@mkboth\markboth
\@ifundefined{chapter}{\def\sectionmark##1{\markboth
{\uppercase{\ifnum \c@secnumdepth>\z@
 \thesection\hskip 1em\relax \fi ##1}}{}}
\def\subsectionmark##1{\markright {\ifnum \c@secnumdepth >\@ne
 \thesubsection\hskip 1em\relax \fi ##1}}}
{\def\chaptermark##1{\markboth {\uppercase{\ifnum \c@secnumdepth>\m@ne
 \@chapapp\ \thechapter. \ \fi ##1}}{}}
\def\sectionmark##1{\markright{\uppercase{\ifnum \c@secnumdepth >\z@
 \thesection. \ \fi ##1}}}}
\def\@oddhead{\@fancyhead\relax\@olhead\@ochead\@orhead\hss}
\def\@oddfoot{\@fancyfoot\relax\@olfoot\@ocfoot\@orfoot\hss}
\def\@evenhead{\@fancyhead\hss\@elhead\@echead\@erhead\relax}
\def\@evenfoot{\@fancyfoot\hss\@elfoot\@ecfoot\@erfoot\relax}
\headwidth\textwidth}
\def\ps@fancyplain{\ps@fancy \let\ps@plain\ps@plain@fancy}
\def\ps@plain@fancy{\@fancyplaintrue\ps@fancy}
\author{J.-P. Serre}\def\authorname{J.-P. Serre}
\newcommand{\R}{\mathbb{R}}
\newcommand{\C}{\mathbb{C}}
\newcommand{\Z}{\mathbb{Z}}
\newcommand{\Q}{\mathbb{Q}}
\newcommand{\F}{\mathbb{F}}
\newcommand{\A}{\mathbb{A}}
\newcommand{\G}{\mathbb{G}}
\renewcommand{\P}{\mathbb{P}}
\renewcommand{\k}{k}

\newcommand{\inprod}[1]{\langle #1\rangle}
\renewcommand{\check}{^{\vee}}
\newcommand{\g}{\mathfrak{g}}
\newcommand{\nilp}{\text{nilp}}
\newcommand{\st}{\ | \ }
\newcommand{\into}{\hookrightarrow}
\newcommand{\twobytwo}[4]{\left(\begin{smallmatrix} #1 & #2 \\ #3 &%
      #4\end{smallmatrix}\right)}
\newcommand{\size}[1]{|#1|}
\newcommand{\urad}{R^{u}}
\newcommand{\sat}{^{\rm sat}}
\renewcommand{\th}{\ensuremath{{}^{\rm th}}}
\newcommand{\tensor}{\otimes}
\renewcommand{\epsilon}{\varepsilon}

\newcommand{\define}[1]{\emph{#1}}
\newcommand{\semidirect}{\rtimes}
\newcommand{\supgroup}{\supset}
\newcommand{\subgroup}{\subset} 
\newcommand{\adjoin}[1]{[#1]} 
\newcommand{\subspace}{\subset}
\newcommand{\SGindex}[2]{(#1:#2)}

\DeclareMathOperator{\Lie}{Lie}
\DeclareMathOperator{\rank}{rank}
\DeclareMathOperator{\SO}{SO}
\DeclareMathOperator{\PGL}{PGL}
\DeclareMathOperator{\PSL}{PSL}
\DeclareMathOperator{\GL}{GL}
\DeclareMathOperator{\SL}{SL}
\DeclareMathOperator{\SU}{SU}
\DeclareMathOperator{\Spin}{Spin}
\DeclareMathOperator{\SP}{Sp}

\DeclareMathOperator{\Aut}{Aut}
\DeclareMathOperator{\Sym}{Sym}
\DeclareMathOperator{\ab}{ab}
\DeclareMathOperator{\Hom}{Hom}
\DeclareMathOperator{\Ext}{Ext}
\DeclareMathOperator{\End}{End}
\DeclareMathOperator{\charistic}{char}
\DeclareMathOperator{\obs}{obs}

\newenvironment{theorem_count}
     {\setlength{\partopsep}{0in}
      \setlength{\itemsep}{0in}
      \setlength{\topsep}{0in}
      \setlength{\parsep}{0in}
      \setlength{\parskip}{0in}
      \renewcommand{\labelenumi}{{\rm (\roman{enumi})}}
      \begin{enumerate}
      \setlength{\parskip}{0in}
      \setlength{\partopsep}{0in}
      \setlength{\itemsep}{0in}
      \setlength{\topsep}{0in}
      \setlength{\parsep}{0in}}
     {\end{enumerate}}
\newcommand{\MakeArabic}
      {\renewcommand{\labelenumi}{{\rm (\arabic{enumi})}}}

\def\dynkinEEE#1#2#3#4#5#6#7#8#9
   {\begin{picture}(124, 50)%
   \put(5, 35){\circle{4}}%
   \put(7, 35){\line(1, 0){15}}%
   \put(24, 35){\circle{4}}%
   \put(26, 35){\line(1, 0){15}}%
   \put(43, 35){\circle{4}}%
   \put(45, 35){\line(1, 0){15}}%
   \put(62, 35){\circle{4}}%
   \put(64, 35){\line(1, 0){15}}%
   \put(81, 35){\circle{4}}%
   \put(83, 35){\line(1, 0){15}}%
   \put(100, 35){\circle{4}}%
   \put(102, 35){\line(1, 0){15}}%
   \put(119, 35){\circle{4}}%
   \put(43, 33){\line(0, -11){15}}%
   \put(43, 16){\circle{4}}%
   \put(5, 41){\makebox(0, 0)[b]{$#1$}}%
   \put(24, 41){\makebox(0, 0)[b]{$#3$}}%
   \put(43, 41){\makebox(0, 0)[b]{$#4$}}%
   \put(62, 41){\makebox(0, 0)[b]{$#5$}}%
   \put(81, 41){\makebox(0, 0)[b]{$#6$}}%
   \put(100, 41){\makebox(0, 0)[b]{$#7$}}%
   \put(119, 41){\makebox(0, 0)[b]{$#8$}}%
   \put(43, 11){\makebox(0, 0)[t]{$#2$}}%
   \put(-60, 35){\makebox(0, 0)[lb]{$#9$}}%
\end{picture}}

\newtheorem*{main_theorem}{Main Theorem}
\newtheorem{theorem}{Theorem}
\newtheorem{property}{Property}
\newtheorem{cor}{Corollary}
\newtheorem{prop}{Proposition}
\theoremstyle{definition}
\newtheorem*{remarks}{Remarks}
\newtheorem*{remark}{Remark}

\begin{document}
\thispagestyle{plain}

\vspace{100mm}
\begin{center}
{\Large University of Oregon Mathematics Department}

\vspace{60mm}

{\Huge\bf Moursund Lectures 1998}

\vspace{30mm}

{\Huge J.-P. Serre}

\vspace{15mm}

\vspace{20mm}

\vspace{20mm}

{Notes by W.E. Duckworth}
\end{center}
\pagebreak
\thispagestyle{empty}
\vspace*{50mm}
These informal
notes are closely based on a series of eight lectures given by J.-P. Serre
at the University of Oregon in October 1998. Professor Serre gave
two talks per week for four weeks. 

The first talk each week
was concerned with
constructing
embeddings of finite groups, especially $\PSL_2(p)$ and $\PGL_2(p)$,
into Lie groups.
The second talk each week was about
generalizations of the
notion of complete reducibility in group theory, especially in positive
characteristic.

The notes are divided into two parts, one for each of the topics
of the lecture series.
At the end of the notes, there is a short list of references
as a guide to further reading.

\pagebreak
\setcounter{page}{1}
\thispagestyle{empty}
\vspace*{20mm}
\begin{center}
{\bf\huge Part I}
\end{center}
\vspace{30mm}
\begin{center}
{\Huge Finite subgroups of Lie groups}
\end{center}
\pagebreak
\vspace{30mm}
\section*{Lecture 1}

We begin with a guiding example. Let $G$ be the compact Lie group
$\SO_3(\R)$. The finite subgroups of $G$ fall into
the following families:
\begin{itemize}
	\item The cyclic subgroup $C_n$ of order $n$. This appears as a
subgroup of the maximal
torus $T$ of $G$ consisting of rotations around some fixed axis. It is not
really interesting: it's there because of the torus, not really because of $G$.
	\item The dihedral group $D_n$ of order $2n$. Again, such subgroups lie
in another Lie subgroup of $G$, namely the normalizer $N$ of $T$ in $G$.
The index $\SGindex{N}{T} = 2$ and the additional reflection generating
$D_n$ lies inside $N$.
	\item Three more ``exceptional'' examples:
the alternating group $A_{4}$ on four letters,
the symmetric group $S_{4}$, and the alternating group $A_{5}$.
These may be viewed as the automorphisms of
the regular tetrahedron, cube and icosahedron respectively.
\end{itemize}

Let us indicate one reason
for the importance of this example for
complex analysis and topology.
One can view $\SO_{3}(\R)$
as a maximal compact subgroup of the group $\PGL_{2}(\C)$, that is,
the group of all
transformations $z \mapsto \frac{az+b}{cz+d}$ with $ad-bc \neq 0$.
Up to conjugacy, compact subgroups of $\PGL_2(\C)$ and $\SO_3(\R)$
are the same. So the above list also describes the embeddings of
finite subgroups $\Gamma$ into $\PGL_2(\C)$.
Now, $\PGL_2(\C)$ is the automorphism group of the projective line
$\P_1$ over $\C$, so a finite subgroup $\Gamma\subgroup \PGL_2(\C)$ acts on
$\P_1$. Dividing, we get a (ramified) Galois covering
$$
\P_1 \stackrel{\Gamma}{\longrightarrow} \P_1 / \Gamma \cong \P_1
$$
of a curve of genus 0 by another, and our list of finite subgroups
gives all possible Galois coverings of $\P_1$ by $\P_1$.

We wish to consider finite subgroups of more general Lie groups $G$.
We will restrict our
attention to the following sorts of Lie group:
\begin{itemize}
    \item Compact, real, connected Lie groups,
especially the semisimple ones: $\SU_n$, $\SO_{n}$, \ldots,$E_{8}$.
    \item The corresponding complex groups: $\SL_{n}(\C)$,
$\SO_{n}(\C)$, \ldots, $E_{8}(\C)$.
    \item Any of these groups $G(\k)$ over an arbitrary field $\k$. Indeed,
thanks to Chevalley, we can define these groups even over $\Z$.
\end{itemize}

In fact as we will see, one can often use the groups $G(\k)$ over
fields of positive characteristic to shed light on the first two
problems.  An example of this philosophy appears in the work of
Minkowski, who was studying lattices $\Lambda\subset \R^{n}$
(cf. [Min]).  The group $\Gamma:=\Aut(\Lambda)$ is finite, and he was
interested in finding an upper bound for the exponent of a given prime
$\ell$ in $\size{\Gamma}$.  Now $\Lambda\cong \Z^{n}$ so
$\Gamma\subgroup \GL_{n}(\Z)$.  If we reduce modulo $p$ then we have a
map $\Gamma\to \GL_{n}(\Z/p\Z)$.  Minkowski showed that for $p\ge3$
this is an embedding, so that $\size{\Gamma}$ divides
$\size{\GL_{n}(\Z/p\Z)} =(p^{n}-1)(p^{n}-p)\dots (p^{n}-p^{n-1})$.
Now, by varying $p$ one gets an upper bound for the exponent of $\ell$
in $\size{\Gamma}$, namely, $\left[ \frac{n}{\ell-1} \right] + \left[
\frac{n}{\ell(\ell-1)} \right] +\ldots$\ .  (This is correct only for
$\ell >2$; the case $\ell=2$ requires a slightly different argument.)
Moreover, this upper bound is exact.\bigskip

From now on $G$ is a semisimple group, e.g. $\SL_{n}$, \ldots,
$E_{8}$.  We want to understand the possible finite groups
$\Gamma\subgroup G(\C)$.  First, we discuss the case that $\Gamma$ is
{\em abelian}.  Let $T$ be a maximal torus of $G$ of dimension $r =
\rank G$.  So, $T\cong \G_{m}\times \ldots \times \G_{m}$ ($r$ copies)
where $\G_m$ is the one dimensional multiplicative group.  So over
$\C$, $T(\C) \cong \C^* \times\dots\times \C^*$.  Thus  we can realize
any abelian finite group on $r$ generators as a subgroup of $T(\C)$.
In fact, \emph{almost} all finite abelian subgroups subgroups of
$G(\C)$ arise in this way, but there are exceptions.  For example,
recall our embedding of the Klein group $D_2$, which is an elementary
abelian $(2,2)$-group, in $\SO_3(\R) \subgroup \PGL_2(\C)$: it
cannot be embedded in $T(\C)$ since $\PGL_2(\C)$ only has rank $1$.

Let us  restrict our attention to elementary abelian $(p, p,
\dots, p)$-groups $E$.  Then all subgroups of $G(\C)$ isomorphic to
$E$ are `toral', that is, are contained in some maximal torus, unless
$p$ is one of finitely many \emph{torsion primes}.  For each of these
torsion primes, there is an `exceptional' embedding of some $E$ into
$G(\C)$; R. Griess has classified such embeddings, cf.  [G].  The torsion
primes for simply connected, simple $G$ are as follows:
\begin{center}
\begin{tabular}{|c|c|c|c|c|c|c|c|c|}
    \hline
    $A_{n}$ &$B_{n}, n\ge 3$ & $C_{n}$ & $D_{n}, n\ge 4$ & $G_{2}$ & $F_{4}$ & $E_{6}$ &
$E_{7}$ &
    $E_{8}$  \\
    \hline
    none & 2 & none & 2 & 2 & 2,3 & 2,3 & 2,3 & 2,3,5  \\
    \hline
\end{tabular}
\end{center}
These primes first arose in topology in the 1950s (cf, e.g. [Bo]). For
a compact Lie group $G$, the cohomology ring $H^*(G, \Z)$ with
coefficients in $\Z$ is not always a free $\Z$-module.  The primes $p$
such that $H^*(G, \Z)$ has $p$-torsion are called the torsion
primes. Moreover, these primes can be described in terms of the root
data: if all the roots have the same length, they are the primes that
divide a coefficient of the highest root when written in terms of the
simple roots.

For another example of one of these exceptional embeddings, consider
$G=G_{2}$.  We view $G_2(\C)$ as the group of automorphisms of the
\emph{Cayley algebra}.  This algebra has a standard
basis $\{1,e_{\alpha}\st \alpha\in \Z/7\Z\}$ and the
automorphisms determined by $e_{\alpha}\mapsto \pm e_{\alpha}$ for
$\alpha=1,2,3$
give an `exceptional' elementary abelian $(2,2,2)$-group inside
$G_2(\C)$.  This subgroup is  important in studying the Galois
cohomology of $G_{2}$.

These  results  on abelian subgroups can be extended
somewhat to \emph{nilpotent subgroups} using a result of Borel-Serre
(cf. [BS]): every finite nilpotent subgroup of $G(\C)$ is contained in
the normalizer of some maximal torus of $G(\C)$.

Now we consider a second, quite different situation, namely, we let
$\Gamma$ be a \emph{quasi-simple} group (i.e. the quotient by its
center is simple and non-abelian).  If $G=\SL_{n}$ then one can
classify all possible embeddings $\Gamma\into G(\C)$ by viewing the
natural $G(\C)$-module as a representation of $\Gamma$ and using
character theory.  A variation of this approach allows one to tackle
the problem also for $G=\SO_{n}$, $\SP_{2n}$ and even $G_2$ (since
this can be viewed as the subgroup of $\SO_{7}$ which leaves invariant
an alternating 3-linear form).  In other words, in these cases, the
problem can be reduced to a question about the character table of
$\Gamma$.  This leaves the cases $F_{4}$, $E_{6}$, $E_{7}$ and
$E_{8}$.  A lot of work in the last few years, in particular by A.
Cohen, R. Griess and A. Ryba, has resulted in a list of the possible
$\Gamma$ that can arise.  This list is complete according to computer
verifications.  There are still open questions however.  For instance,
the number of conjugacy classes of such subgroups is not known in
general.

Some of the most interesting questions arise when
$\Gamma=\PSL_{2}(\F_{p})$.  For instance, if $G= E_8$, then $G(\C)$
has finite subgroups $\PSL_2(\F_{p})$ for $p=31,41,61$ (cf.[GR], [S3]).  The
principal difficulty is in proving the \emph{existence} of these subgroups.
We now discuss briefly the sorts of method one can use for such a
construction.

The first method  depends upon  computer
calculations.  For instance, to embed $\PSL_{2}(\F_{61})$ in
$E_{8}(\C)$, start with a Borel subgroup $B \subgroup \PSL_{2}(\F_{61})$
consisting of all upper triangular matrices and its opposite $B^{-}$
consisting of all lower triangular matrices.  Then $B$ is isomorphic
to a semidirect product of cyclic subgroups of order 61 and 30.  Choose
an element in $E_{8}(\C)$ of order 30, namely, a \emph{Coxeter
element}.  There is a subgroup of $E_8(\C)$ generated by
an element of order 61 upon which this Coxeter element acts by an
automorphism of order 30.  We map
$B$ to the subgroup of $E_8(\C)$ generated by these two elements.
Then one needs an involution within $E_8(\C)$ which gives the
embedding of the other Borel $B^-$, and this is where the computer
comes in.  In fact, the computer calculations are done by working
within $E_{8}(\F_{\ell})$ for some large prime $\ell$ not dividing the
order of $\Gamma$. The results are then lifted
(easily) to $E_{8}(\C)$.

The second method (cf.[S3]) is quite different, and depends on lifting from
 \emph{the same characteristic} $p = 61$.  There is a so-called
\define{principal homomorphism}
$\SL_{2}\rightarrow E_{8}$
with kernel $\{\pm 1\}$.  This is defined over $\F_p$, giving
an embedding of $\PSL_2(\F_p)$ into $E_8(\F_p)$.
The idea is to lift this embedding to
an embedding in characteristic $0$.  However, there may be a non-trivial
 obstruction preventing a lift to an embedding
 $\PSL_{2}(\F_{p})\into E_{8}(\Z / p^2 \Z)$.  So one has to
 proceed more indirectly, and we will discuss the argument in more
 detail in the remaining lectures.

Finally, we return to our opening example. Recall we had
subgroups $A_{4} \cong \PSL_{2}(\F_{3})$, $S_{4} \cong
\PGL_{2}(\F_{3})$ and
$A_{5}\cong \PSL_{2}(\F_{5})$
inside $\SO_{3}(\R)$, corresponding to the symmetries
of the tetrahedron, cube and icosahedron.
The analogues of these embeddings for
$E_{8}$ are the embeddings of $\PSL_{2}(\F_{31})$, $\PGL_{2}(\F_{31})$
and $\PSL_{2}(\F_{61})$! In fact, quite generally for any simple,
simply-connected $G$, let
$h$ be the \define{Coxeter number} defined as
$\dfrac{\dim G}{\rank G} -1$.  Notice in the rank 1
case $\SO_{3}(\R)$, we have $h=2$.  In the case of $E_{8}$ we
have $h=30$.  It is true in general that if $h+1$ or $2h+1$
is prime then $G(\C)$ has subgroups of
the form
$\PSL_{2}(\F_{h+1})$,
$\PGL_{2}(\F_{h+1})$ and
$\PSL_{2}(\F_{2h+1})$.

\section*{Lecture 2}
We continue to assume that $G$ is a simple algebraic group
over an algebraically closed field $\k$ of characteristic zero.
We recall our notation:  $r = \rank G$,  $h$ is the Coxeter
number $\frac{\dim G}{r}-1$, and $W_{G}$ is the Weyl group of $G$ (uniquely
determined up to isomorphism).
Also fix $q = p^e$ for some prime $p$.

The group $W_G$ has a natural reflection representation $V$ of
dimension $r$.  Let $\k[V]$ denote the coordinate ring of $V$, a
polynomial ring in $r$ generators.  By general theory (cf.[B], Chap
V, \S5), the ring of invariants $\k[V]^{W_G}$ is a graded
polynomial ring in $r$ generators, $P_1, \dots, P_r$ say.
Moreover, the degrees $2 = d_1 \leq d_2 \leq \dots \leq d_r = h$ of
these generators $P_1, \dots, P_r$ are uniquely determined.  
The invariant degrees are listed in Table 1.
\begin{table}
\begin{center}
{\small\sc Table 1: Invariant degrees}
\end{center}\begin{center}
\begin{tabular}{|c|c|c|c|}
    \hline
    $G$ & degrees & $\dim G$ & $h$\\
    \hline
    $A_{r}$ & 2,3,\ldots,$r+1$ & $(r+1)^{2}-1$ & $r+1$  \\
    $B_{r}$ & 2,4,\ldots,$2r$ & $2r^{2}+r$ & $2r$  \\
    $C_{r}$ & 2,4,\ldots,$2r$ & $2r^{2}+r$ & $2r$  \\
    $D_{r}$  & 2,4,\ldots,$2r-2$,$r$ & $2r^2 - r$ & $2r-2$  \\
    $G_{2}$ & 2,6 & 14 & 6  \\
    $F_{4}$ & 2,6,8,12 & 52 & 12  \\
    $E_{6}$ & 2,5,6,8,9,12 & 78 & 12  \\
    $E_{7}$ & 2,6,8,10,12,14,18 & 133 & 18  \\
    $E_{8}$ & 2,8,12,14,18,20,24,30 & 248 & 30  \\
    \hline
\end{tabular}\end{center}
\end{table}

Now let $\Gamma$ be
either $\SL_{2}(q)$ or $\GL_{2}(q)$.  Let $U$ be the unipotent
subgroup of $\Gamma$ consisting of all upper triangular unipotent matrices,
so $U$ is an elementary abelian group of type
$(p,\ldots,p)$ ($e$ times).  Suppose we have a map
$$
f :\Gamma \rightarrow G,
$$
which is {\em nondegenerate} in the sense that
$\ker f$ is contained in the center of $\Gamma$.
We say that $f$ is of \define{toral type} if $f(U)$ is contained
in a torus of $G$.

In the remaining lectures we will give a partial proof of the following:

\begin{main_theorem} Suppose $q\ge 5$.  There exists
a nondegenerate map $$f:\SL_{2}(q)\to G$$ of toral type
if and only if $q-1$ divides $2d$ for some degree $d$.
\end{main_theorem}

We begin with the easy implication, namely, that
the existence of such a map $f$ implies that $q-1$ divides some $2d$.
In fact one proves more: if there is a nondegenerate toral
map from  a Borel subgroup of $\Gamma = \SL_2(q)$ to $G$ then
$q-1$ divides $2d$.

Let $T$ be a maximal torus in $G$, $N$ its normalizer.  Then
$N/T=W_{G}$ acts on $T$.  Moreover, $N$ \emph{controls
the fusion} of $T$ in $G$; this means:

\vspace{2mm}

\noindent{\bf (F) . }{\em
If $A$ and $A'$ are subsets of $T$, $g\in G$ with
$gAg^{-1}=A'$ then there exists $n\in N$ such that $nan^{-1}=gag^{-1}$
for all $a\in A$.}

\vspace{2mm}

We will also need the following theorem of Springer [Sp2]:

\vspace{2mm}

\noindent{\bf (Sp). }
{\em Let $m\ge 1$.  The following are equivalent:
\begin{theorem_count}
\item   $m$ divides one of the degrees of $W_{G}$;

\item there exists $w\in W_{G}$ and an eigenvalue $\lambda$ of $w$
    (for the natural representation) whose order is $m$.
\end{theorem_count}}

We assume now that $q$ is odd (the even case being similar).
Let $B$ be the Borel subgroup of
$\SL_{2}(q)$ consisting of all upper triangular matrices.
Let $f$ be a nondegenerate homomorphism of
$B$ into $G$.  Let $A$ be the image of $U$; we
can assume that $A\subgroup T$.
Then
$A\cong \F_{q}$ and $B$ acts upon $A$ as squares, that is,
conjugating by $\twobytwo t00{t^{-1}}$ acts as multiplication by
$t^{2}$.  In particular there exists an automorphism $\sigma$ of $A$
induced by $B$ of the form $a\mapsto \lambda a$ for $\lambda\in
\F_{q}^{*}$ of order $\frac{q-1}{2}$.

Now, $\sigma$ has an
eigenvalue in $\overline{\F}_{p}$ of order $\frac{q-1}{2}$,
and by (F) above,
the action of $\sigma$ is induced by an
element $w\in W_{G}$.  Viewing $A$ as a subset of $T[p]$, which is
the reduction modulo $p$ of the standard representation of $W_{G}$, we
deduce that in characteristic zero there exists an eigenvalue of $w$
which has order $\frac{q-1}{2}p^{\alpha}$ for some $\alpha$.  By
(Sp), $\frac{q-1}{2}p^{\alpha}$ divides some degree $d$, as
required.

Note that
the same arguments apply (with minor modifications)
to maps from $\Gamma = \GL_2(q)$ to $G$ as well: in this case,
one finds that $q-1$ divides one of the degrees.

\vspace{2mm}

Now we turn to the converse.
The case where $G$ is classical can be handled directly using the knowledge
of the character table of $\Gamma = \PSL_2(\F_q)$. For instance, if
$G$ is of type $A$, one uses the irreducible representation of $\Gamma$ of degree
$\frac{q-1}{2}$ (assuming $p \neq 2$; if $p = 2$, use a representation of
degree $q - 1$).

So let $G$ be exceptional. One can easily work out
which $\PSL_2(q)$ need to be constructed,
remembering our assumption $q \geq 5$:
\begin{itemize}
    \setlength{\itemsep}{0in}
    \item
    For $G_{2}$, $q-1$ should divide 4 or 12.
    But there is a subgroup $A_2$ of $G_2$, and the case $q-1$ divides
    6 has been treated already, working inside this $A_2$.
    So one just needs to embed
    $\PSL_2(13)$ into $G_2(\k)$.
    \item
    Again, for $F_4$, $q-1$ should divide $4,12,16$ or $24$, but most
    of the first
    two cases have already been dealt with since
    $F_4$ contains a subgroup $G_2$.
    So we need  embeddings of
    $\PSL_2(17)$ and $\PSL_2(25)$.
    \item
    For $E_{6}$ the new cases
    are $q = 11, 19$.
    \item
    For $E_{7}$, they are $q = 29, 37$.
    \item
    For $E_{8}$, they are $q = 31, 41, 49, 61$.
\end{itemize}

We will give a uniform proof of existence in all these cases provided
$q$ is prime. The missing cases (essentially, $q = 25$ for $F_4$ and
$49$ for $E_8$) have been done by computer calculation, cf. [GR].

Some can be done right away with the next theorem, which
for instance covers $E_8$ for $q = 31$.

\begin{theorem} {\rm ([S3])} Let $\k$ be an algebraically closed field of
characteristic 0.  If $p=h+1$ is prime, then there exists a nondegenerate
toral map $\PGL_{2}(\F_{p})\to G$.
\label{exists:map}
\end{theorem}

Let us sketch the proof.  We may assume that $G$ is split,
so that $G(R)$ makes sense for any ring $R$.  In particular we have
$G(\Z/p\Z)$, $G(\Z/p^{2}\Z)$, etc...  and
$$
\varprojlim
G(\Z / p^r\Z) = G(\Z_p)
$$
where $\Z_p$ denotes the ring of $p$-adic integers.  Let us start from
an embedding of $\PGL_2(\F_p)$ into $G(\F_p)$ in which the non-trivial
elements of $U$ are regular unipotent elements of $G(\F_p)$.  The
existence of such an embedding over $\F_p$ was proved by Testerman
(this requires $p \geq h$ which is true in our setting: $p = h+1$)
see [Te] and [S3].

This embedding
of $\Gamma$ into $G(\F_p)$ lifts to $G(\Z / p^2 \Z)$:
$$\begin{matrix}
    & 1 & \to &\Lie G_{/\F_{p}} &\to &G(\Z/p^{2}\Z) &\to &G(\Z/
p\Z)&\to& 1\\
&&&&&&&\uparrow\\
&&&&&&& \Gamma
\end{matrix}$$
Indeed, the obstruction to such a lift is 0 because of:
\begin{theorem}
    $H^{i}(\Gamma,\Lie G_{/\F_{p}})=0$ for $i\ge 1$.
    \label{cohom:zero}
    \end{theorem}
The  proof of Theorem \ref{cohom:zero} uses the embedding
$$
H^{i}(\Gamma,\Lie G_{/ \F_p})\into H^{i}(C_{p},\Lie G_{/\F_{p}}),
$$
where $C_p \cong U$ is a Sylow $p$-subgroup of $\Gamma$.  Now,
$$
\dim
H^{0}(C_{p},\Lie G_{/\F_p}) = \dim \hbox{(Lie algebra of the
centralizer of
$C_{p}$)}$$
and, since the non-trivial elements of $C_p$ are regular, this
dimension is $r$, cf. [St].  Using the fact that $\dim G=pr$, one sees that
every Jordan block of the action of $C_p$ on $\Lie G_{/\F_p}$
has size $p$, and $H^{i}(C_{p},\Lie  G_{/\F_{p}})=0$ as required.

Hence, the lifting to $\Z / p^2 \Z$ is possible. The same argument
applies to $\Z / p^3 \Z$, etc. One ends up with an embedding of
$\PGL_2(\F_p)$ in $G(\Z_p)$, hence in $G(\Q_p)$. Since $\Q_p$
is of characteristic 0, an easy argument then gives an embedding in
$G(\C)$ (or even in $G(K)$ where $K$ is a number field),
as was to be shown. 

\begin{remark} In Theorem \ref{exists:map}, the hypothesis that $\k$ has
characteristic 0 can be suppressed (cf.[S3]), except in one case: $G$
of type $A_{1}$, and $\k$ of characteristic equal to 2.  (Indeed,
there is no embedding of $S_4$ into $\PGL_2(\k)$ when $\charistic (\k)=2$.)
\end{remark} 
\section*{Lecture 3}

Let us  give a  sketch of an existence
proof in the
remaining cases of the Main Theorem with $q$ prime,
postponing some of the technical details until Lecture 4.

As before let $G$ be quasi-simple and  split over $\Z$.  Let $h$ be
its Coxeter number.
Suppose we have a non-trivial morphism $\phi:\SL_{2}\to G$ such that:
\begin{theorem_count}
    \MakeArabic
    \item $\phi$ is defined over the local ring of $\Z$ at $p$,
    i.e. over $\Z_{(p)}$;
    \item  writing  $\Lie G=\bigoplus L(n_{i})$ where $L(n_{i})$ is the
    irreducible representation of $\SL_{2}$ with highest weight
    $n_{i}$,  we require that all $n_{i}$ are $<p$,  exactly one $n_{i}$
equals 2, and exactly one $n_{i}$ equals
    $p-3$;
    \item  $p>h$.
\end{theorem_count}
We will prove:

\begin{theorem} If $\phi:\SL_{2}\to G$ is a morphism satisfying the above
conditions, then there exists a non-degenerate morphism
$\SL_{2}(\F_{p})\to G(\C)$.
\label{exists:map:fromSL2}
\end{theorem}

As a special case take $\phi$ to be the principal embedding, as
discussed by Kostant and others [K].  Here property (1) has been
verified by Testerman [Te].  Denoting the invariant degrees $d_1, \dots,
d_r$, the $n_i$ in this case are $$\{2 d_i - 2 \st i = 1, \dots, r\}.$$
In particular the largest $n_i$ is $2h-2=p-3$ and the conditions (2)
and (3) are satisfied.  We then obtain as a consequence of Theorem
\ref{exists:map:fromSL2} a
proof of a well known
conjecture of Kostant (in the special case where $2h+1$ is prime).

As another special case  consider $G=\SL_{n}$, with degrees are
$2,3,\ldots,n$.  The $n_{i}$'s are 
$\{2d_{i}-2 \st i = 1, \dots,
r\}$.  So taking $n=\frac{p-1}{2}$, the conditions of the theorem are
satisfied, and we recover the well known fact (due to Frobenius)
that $\SL_{2}(\F_{p})$ has an irreducible character of degree
$\frac{p-1}{2}$.  The existence of an irreducible character of
degree $\frac{p+1}{2}$ can be proved similarly.

Other examples come from Dynkin's classification of $A_1$ type
subgroups of simple algebraic groups in characteristic $0$ (see [Dy]).
Dynkin's work shows that such embeddings are determined uniquely up to
conjugacy in the following way.  Let $\phi:\SL_{2}\into G$ be an
embedding and $\{\alpha_{1},\ldots,\alpha_{r}\}$ a base of the root
system of $G$.  We may assume that $\phi$ maps the maximal torus
$\G_{m}$ of $\SL_{2}$ into the maximal torus $T$ of $G$.  For a root
$\alpha$,  the inner product $\inprod{\phi,\alpha}$ is defined as
the integer corresponding to the composite function
$$\G_{m}\stackrel{\phi|_{\G_{m}}}{\longrightarrow}
T\stackrel{\alpha}{\longrightarrow} \G_{m}.$$  
We may also assume that $\phi$ belongs to the Weyl chamber, i.e. that
all $\inprod{\phi,\alpha_i}$ are $\ge 0$.  Then the embedding $\phi$ is
determined up to conjugacy by the weights $\inprod{\phi,\alpha_{i}}$
for $i = 1, \dots, r$.  Writing these on the corresponding nodes of
the Dynkin diagram of $G$, we obtain a \define{labelled diagram}
determining the embedding $\phi$.  Dynkin worked out precisely which
labelled diagrams can arise.  We mention two examples with
$G = E_8$ when the  labelled diagrams are:
$$
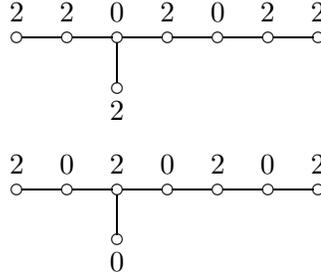
\dynkinEEE{2}{2}{2}{0}{2}{0}{2}{2}
$$
$$
\dynkinEEE{2}{0}{0}{2}{0}{2}{0}{2}
$$
One shows that Theorem \ref{exists:map:fromSL2}, applied to such
diagrams,  gives embeddings with $p = 41$
and $p = 31$.
Similarly, one gets $p=29$ and $p=37$ for $E_{7}$.  (All these cases
have also been done by computer, except $p=29$.)

We now begin the proof of the theorem.  Let $\Q_{p}$ be the field
of $p$-adic numbers.  It is not possible to work over $\Q_{p}$ as, for
example, the values of the character
of $\SL_{2}(p)$ of degree $\frac{p-1}{2}$ involve
$\frac{-1+ \sqrt{\pm p}}{2}$.  So we need to work over the ramified extension
$K_{p,u}:=\Q_{p}(\sqrt{pu\,})$ where $u$ is a unit in $\Z_{p}$
(there are only
two cases according as $u$ is square $\operatorname{mod} p$ or not).
Set $R_{p,u}:=\Z_{p}[\sqrt{pu\,}]$, the corresponding ring of
integers, with residue field $\F_{p}$ as before.
We will  prove:

\begin{theorem} One may choose $u$ so that the subgroup
$\phi(\SL_{2}(\F_{p})) \subgroup  G(\F_{p})$ can be lifted to
a subgroup of $G(R_{p,u})$.
\label{choose:u:to:lift}
\end{theorem}

Viewing $G(R_{p,u})$ as a subgroup of $G(\C)$ this implies Theorem
\ref{exists:map:fromSL2}.

To prove Theorem \ref{choose:u:to:lift}, we first abbreviate
$R=R_{p,u}$, $\pi=\sqrt{pu}$,
$A=\phi(\SL_{2}(\F_{p}))$.  As $G$ is smooth, we have surjective
maps $G(R)\to
G(R/\pi^{n}R)$ with kernels denoted $G_n$.
Then
$G=G_{0}\supgroup G_{1}\supgroup \ldots$ and
$G=\varprojlim G/G_{n}$.
We note the following basic properties (cf. [DG]):
\begin{align*}
    G/G_{1} & =  G(\F_{p})\\
    G/G_{2} & =  G(R/\pi^{2}R)=G(R/pR) \cong  \Lie_{p}G\semidirect G(\F_{p})\\
    (G_{i},G_{j}) & \subgroup G_{i+j}\\
    G_{i}/G_{i+1} & \cong  \Lie_{p}G
    \end{align*}
By assumption, $A$ is embedded in $G / G_1$, and we would like to lift
this to $G/G_{2}$.  The split exact sequence
$$1 \to \Lie_{p} G \to G/G_{2} \to
G/G_{1}\to 1$$
gives an obvious lift
$\sigma:A\to G/ G_2$.
However, this $\sigma$
does not  lift to $G / G_3$,  so we need to modify it.
For any $\alpha \in
H^{1}(A,\Lie_{p}G)$ represented by a 1-cocycle $a$, we can define
a new lift $\sigma_{a}(s)$ by setting $\sigma_a(s) =
\sigma(s)a(s)$.  Two liftings are conjugate
by an element in the kernel if and only if the corresponding cocycles are cohomologous. So
studying $H^{1}(A,\Lie_{p}G)$ is crucial.  We will show that there
is a choice of $\alpha$ which allow us to continue  lifting to every
$G/G_{n}$.

Hypothesis (2) is known to imply that $\Lie_{p}G=\bigoplus L(n_{i})_{p}$ with $n_{i}<p$.
When $n<p$ one can show that $\dim H^{1}(A,L(n))$ is 1 if $n=p-3$ and
zero otherwise.  When $n<p-1$ one can show that $\dim H^{2}(A,L(n))$
is 1 if $n=2$ and zero otherwise.  Hence 
$$\dim H^{1}(A, \Lie_{p}G)=1\quad\text{and}\quad\dim H^{2}(A,\Lie_{p}G)=1.$$
We wish to use the sequence $\Lie_{p}G \to G/G_{3} \to G/G_{2}$ to
lift $A \stackrel{\sigma}{\to}G/G_{2}$ to a map $A\to G/G_{3}$.
The corresponding obstruction is denoted by
$\obs(\sigma)\in H^{2}(A,\Lie_{p}G)$.  With $\sigma$ equal to the lift
of the original embedding $A\into G/G_{1}$ we have
$\obs(\sigma)\ne 0$.
Now if we take $\alpha\in H^{1}(A,\Lie _{p}G)$ and calculate
$\sigma_{\alpha}$ one finds (cf. Lecture 4) that  $\obs (\sigma_{\alpha})
=\obs(\sigma)+\frac{1}{2}[\alpha,\alpha]$ where $[\ ,\ ]$ is the cup
product in cohomology induced by 
$$[\ ,\ ]: \Lie_{p}G\times \Lie_{p}G
\to \Lie_{p}G.$$  
So we try to choose $\alpha$ such that
$\obs(\sigma)+\frac{1}{2}[\alpha,\alpha]=0$.   If there is no
$\alpha$ satisfying this equation, then we change our choice of $u$
(in fact in all the cases I know, the choice of $u=-1$ works).

So now we may assume that $A\stackrel{\sigma}{\to}G/G_{2}$ is liftable to
$G/G_{3}$.  Call this lift $\tau$.  The lift to $G/G_{4}$  may
again have a non-trivial  obstruction in $H^{2}(A,\Lie_{p}G)$.  Again
one can modify $\tau$ by a 1-cocycle $b:A\to \Lie_p G$, and one proves
that $\obs(\tau_b)=\obs (\tau) +[\alpha, \beta]$, where $\beta$ is the
class of $b$ in $H^1(A,\Lie_pG)$.  Since $\alpha\ne 0$ one can choose
$\beta$ such that $[\alpha,\beta]=-\obs (\tau)$, hence
$\obs(\tau_b)=0$ and $\tau_b$ can be lifted to $G/G_4$.

The process continues in this manner:
we obtain inductively a lift of a map
$A\to G/G_{n}$ to $G/G_{n+1}$ and then modify this lift to
get a map to $G/G_{n+2}$.
Putting it all together completes our sketch of the
proof of Theorem \ref{choose:u:to:lift}.

\section*{Lecture 4}

In this lecture, we go back to discuss some of the  technical points
arising in the proof of Theorem \ref{choose:u:to:lift}.
We will consider a general setup which includes the situation considered
in Lecture 3.

Consider a sequence of  surjective group homomorphisms $E_{3}\to E_{2}\to
E_{1}$ with $M_{1}=\ker E_{3}\to E_{2}$, $M_{2}=\ker E_{3}\to E_{1}$
and $M_{3}=\ker E_{2}\to E_{1}$.  One has a short exact sequence:
$$1\to M_{1}\to M_{2}\to M_{3}\to 1.$$

\vspace{1mm}
\noindent{\em Assumption A.} Assume that $M_{1}$, $M_{3}$ are abelian,  and
$M_{1}$ is in the center of $M_{2}$.  (This gives natural actions of
$E_{2}$ on $M_{2}$ and of $E_{1}$ on $M_{1}$ and $M_{3}$.)
\vspace{1mm}

Now let $A$ be a group and $\phi:A\to E_{2}$.  Call $\obs(\phi)\in
H^{2}(A,M_{1})$ the obstruction to lifting $\phi$ to $A\to E_{3}$.
Let $x$ be a 1-cocycle $A\to M_{3}$, and $\phi_{x}$ be the map $s\mapsto
x(s)\phi(s)$ of $A$ into $E_2$. Write $\underline{x}$ for the class of $x$ in
$H^{1}(A,M_{3})$.  We want to compare $\obs(\phi)$ and
$\obs(\phi_{x})$; note that $A$ acts the same way on $M_{1}$ by $\phi$
or by $\phi_{x}$ since $E_{1}$ acts on $M_{1}$.
We have the following key formula:
\begin{prop}
     $\obs(\phi_{x}) = \obs(\phi)+\Delta(\underline{x})$
     \label{formula:for:obs}
\end{prop}
\noindent
where $\Delta: H^{1}(A,M_{3})\to H^{2}(A,M_{1})$ is the (non-abelian)
coboundary map associated with the exact sequence of $A$-groups:
$$1\to M_{1}\to M_{2} \to M_{3} \to 1,$$
(cf. [S4] Ch I, \S5.7).
This formula will be verified by
a direct computation given at the end of the lecture.

Next, we want to compute $\Delta: H^{1}(A, M_{3}) \to H^{2} (A,
M_{1})$.  We make the following assumption:

\vspace{1mm}
\noindent{\em Assumption B.}
The map  $m\mapsto m^{2}$ of $M_{1}$ onto
itself is bijective (in additive notation $M_{1}$ is a $\Z[\frac
12]$-module).
\vspace{1mm}

This allows us to define an \emph{addition} in $M_{2}$ by
$x+y=x.y.(x,y)^{-1/2}$ (note that $(x,y)$ is the usual commutator and
belongs to $M_{1}$) and a {\em Lie bracket} $[x,y]=(x,y)$.
This makes $M_{2}$
into \emph{a Lie algebra}.  (I am using here an elementary case of
the inversion of the Hausdorff formula, cf. [B], Chap II, \S6.)

Call $M_{2}^{\ab}$ the corresponding abelian group.  We have an exact
sequence of $A$-modules:
$$
0\to M_{1}\to M_{2}^{\ab}\to M_{3} \to 0,
$$
hence by (abelian) cohomology an additive map
$$
\delta:H^{1}(A,M_{3})\to H^{2}(A, M_{1}).
$$
On the other hand the
bracket defines a bilinear map
$$
M_{3}\times M_{3}\to M_{1},
$$ hence a
(cup product) map
$$
H^{1}(A, M_{3})\times H^{1}(A, M_{3})\to H^{2}(A,M_{1})
$$
which we denote by $\alpha,\beta\mapsto [\alpha.\beta]$.  It is
symmetric.  We can now state the formula giving $\Delta$:
\begin{prop}
    $\Delta(\alpha)=\delta(\alpha)+\frac12[\alpha.\alpha]$
    \label{formula:for:Delta}
for every $\alpha\in H^{1}(A, M_{3})$.
\end{prop}  This is verified again by a
direct computation which we will give at the end of the lecture.
Note that $\delta(\alpha)$ is linear in $\alpha$ and $[\alpha.\alpha]$
is quadratic; hence $\Delta$ is a polynomial function of degree 2.
We conclude with the computations mentioned above.

\vspace{4mm}

\noindent
\emph{Computations for Proposition \ref{formula:for:obs}.}
If $s\in A$, one has $\phi(s)\in E_{2}$; choose $z_{s}\in E_{3}$,
with $z_{s}\mapsto \phi(s)$.  This defines a 2-cocycle $o(s,t)$ by the
usual formula
$$
z_{s}z_{t}=o(s,t)z_{st}\ , \quad o(s,t)\in M_{1}.
$$
The class of $o(s,t)$ in $H^{2}(A, M_{1})$ is $\obs(\phi)$.

Similarly, choose $b_{s}\in M_{2}$ with $b_{s} \mapsto x(s)$ in
$M_{1}$.  By [S4], \emph{loc. cit.}, $\Delta(\underline{x})$ is the class
in $H^{2}(A,M_{1})$ of the 2-cocycle $\Delta(s,t)$ defined by
$$
\Delta(s,t)=b_{s}.{}^{s}b_{t}.b_{st}^{-1}
$$
(where ${}^{s}b_{t}$ means the transform of $b_{t}$ by $\phi(s)$, i.e.
$z_{s}b_{t}z_{s}^{-1}$.)  Since
$\phi_{x}(s)=x(s)\phi(s)$ we may choose $b_{s}z_{s}$ as a lifting
of $\phi_{x}(s)$ in $E_{3}$.  This gives a cocycle $o_{x}(s,t)$ by:
$$
b_{s}z_{s}.b_{t}z_{t}=o_{x}(s,t).b_{st}z_{st}
$$
and the class of $o_{x}(s,t)$ is $\obs(\phi_{x})$.

We calculate $b_{s}z_{s}b_{t}z_{t}$:
\begin{eqnarray*}
    b_{s}.z_{s}b_{t}z_{s}^{-1}.z_{s}z_{t} & = &
b_{s}.z_{s}b_{t}z_{s}^{-1}.o(s,t)z_{st} \\
& = &
b_{s}.{}^{s}b_{t}.o(s,t).z_{st} \\
&=& \Delta(s,t)b_{st}o(s,t)z_{st}.
\end{eqnarray*}
Hence $o_{x}(s,t) b_{st} = \Delta(s,t) b_{st}o(s,t)$.  Since $b_{st}$
commutes with $o(s,t)$, this gives $o_{x}(s,t) = \Delta(s,t).o(s,t)$, as
desired.
\vspace{4mm}

\noindent
\emph{Computations for Proposition \ref{formula:for:Delta}.}
Choose a 1-cocycle $(a_{s})$ of $A$ in $M_{3}$ representing the
class $\alpha$, and lift $a_{s}$ to $b_{s}\in M_{2}$.  The cocycle
$\Delta(s,t)$ defined by
$$
\Delta(s,t)= b_{s}.{}^{s}b_{t}.b_{st}^{-1}
$$
represents $\Delta(\alpha)$, cf. above.

On the other hand, the coboundary $\delta(\alpha)$ may be represented
by the 2-cocycle $\delta(s,t)$ given by
$$
\delta(s,t) = b_{s} * {}^{s}b_{t} * b_{st}^{-1},
$$ where $x*y$ is the product of $x,y$ with respect to the
composition law $x.y.(x,y)^{-1/2}$.  By collecting terms, this gives
$$
\delta(s,t)= \Delta (s,t) \gamma(s,t),
$$
where $\gamma(s,t)=(b_{s},{}^{s}b_{t})^{-1/2} (b_{s}{}^{s}b_{t},
b_{st}^{-1})^{-1/2}$.  In additive notation, this means:
$$
\gamma(s,t)=-\frac12
[a_{s},{}^{s}a_{t}]+\frac12[a_{s}+{}^{s}a_{t},a_{st}].
$$
But $a_{s}+{}^{s}a_{t}=a_{st}$, since $a$ is a 1-cocycle.  Hence the
last term is 0.  As for $s,t\mapsto [a_{s}, {}^{s}a_{t}]$, it is the
cup-product (with respect to $[\ ,\ ]$) of the cocycle $a$ with
itself.  Hence $\gamma=-\frac12[\alpha.\alpha]$ and since
$\Delta(\alpha) =\delta(\alpha)-\underline{\gamma}(\alpha)$, where
$\underline{\gamma}(\alpha)$ is the class of $\gamma(s,t)$, this gives
the required formula.

\pagebreak
\thispagestyle{empty}
\vspace*{20mm}
\begin{center}
{\bf\huge Part II}
\end{center}
\vspace{30mm}
\begin{center}
{\Huge The notion of complete reducibility in group theory}
\end{center}
\pagebreak
\setcounter{theorem}{0}
\setcounter{cor}{0}
\vspace{30mm}
\section*{Lecture 1}

Let $\Gamma$ be a  group.  We will discuss \emph{linear
representations}
of $\Gamma$ over some fixed field $\k$ of characteristic $p \geq 0$.
By this we mean a group
homomorphism $\Gamma\to \GL(V)$ for some finite dimensional vector
space $V$ over $\k$.  We will usually refer to $V$ instead
as a $\Gamma$-module,
though of course technically we should say
 $\k \adjoin\Gamma$-module where $\k \adjoin\Gamma$ denotes the group
algebra of $\Gamma$
over $\k$.  Recall that $V$ is
\emph{irreducible} or \emph{simple} if:
\begin{theorem_count}
    \MakeArabic
    \item  $V\ne 0$;
    \item no subspace of $V$ is $\Gamma$-stable apart from 0 and $V$.
\end{theorem_count}
One says that $V$ is \emph{completely reducible} or \emph{semisimple}
if $V$ is a direct sum of irreducible submodules; equivalently, $V$ is
semisimple if $V$ is generated by irreducible submodules.

The category of semisimple $\Gamma$-modules is stable under the usual
operations of linear algebra.  In other words one can take
$\Gamma$-stable subspaces, quotients, direct sums and duals all within
this category.  Indeed, all of these statements (apart from dual
spaces) are true for modules over an arbitrary ring.  But when we
consider groups, we can also consider the operations of multilinear
algebra.  For instance, given two $\Gamma$-modules $V_{1},V_{2}$ we
can impose a $\Gamma$-module structure upon $V_{1}\tensor V_{2}$ using
the diagonal map $\Gamma \to \Gamma \times \Gamma$.  From this
we can construct exterior powers, symmetric powers, etc....

Around 1950,  Chevalley  proved the following simple looking result:
\begin{theorem} {\rm (cf. [C])} Suppose that $\k$ has characteristic $0$.
If $V_{1}$,
$V_{2}$ are
semisimple $\Gamma$-modules, then $V_{1}\tensor V_{2}$ is again
semisimple.
\label{chevalleys:theorem}
\end{theorem}

An interesting feature of this result is that, although it is stated
in elementary terms, the only known proofs involve some algebraic
geometry.  We sketch the idea.  One starts with a series of
reductions, reducing to the case that $\k$ is algebraically closed and
$\Gamma$ is a subgroup of $\GL(V_1)\times\GL(V_2)$.  Then one replaces
$\Gamma$ by its Zariski closure in $\GL(V_1)\times \GL(V_2)$.  So now
$\Gamma$ is an algebraic group.  The connected component
$\Gamma^\circ$ of $\Gamma$ containing the identity is a normal
subgroup of $\Gamma$ of finite index (this is one bonus of using the
Zariski topology).  In other words, $\Gamma / \Gamma^\circ$ is a
finite group and since the characteristic is 0, one easily then
reduces to the case that $\Gamma = \Gamma^\circ$.  So now, $\Gamma$ is
connected.  Let $\urad\Gamma$ be the unipotent radical of $\Gamma$,
i.e. its  largest normal unipotent subgroup.  In any semisimple
representation, $\urad\Gamma$ acts trivially, and the converse is
known to be true in characteristic zero.  Since $V_1$ and $V_2$ are
semisimple and the representation of $\Gamma$ on $V_1 \oplus V_2$ is
faithful, we deduce that $\urad \Gamma$ is trivial, and we are done.

Now we ask what happens for $p > 0$.  Chevalley's result does not
remain true in general.  For instance, consider
$\Gamma=\SL_{2}(\k)=\SL(V)$ with $\dim V=2$.  Let $\Sym^{n}(V)$ be  $n\th$
symmetric power of
$V$.  If $n<p$ then $\Sym^{n}(V)$ is an irreducible
representation of $\Gamma$.  But if $n=p$, the subspace
$V^{[1]}\subspace \Sym^{p}(V)$ generated by $x^{p}$ and $y^{p}$, where
$\{x,y\}$ is any basis of $V$, is stable under the action of $\Gamma$.
This gives a short exact sequence $$0\to V^{[1]} \to \Sym^{p}(V)\to
L\tensor\Sym^{p-2}(V) \to 0$$
where $L=\det V$ is one-dimensional.  This sequence does
not split (unless $p=\size{\k} = 2$).
So $\Sym^{p}(V)$ is not semisimple in general. Hence,
$V\tensor\ldots\tensor V$ ($p$ times) is not semisimple either.

Now a general principle is that if a
statement is true in characteristic zero then it is also true for
``large'' $p$.  In keeping with this, we have the following:
\begin{theorem} {\rm ([S1])} Let $V_{1},\ldots,V_{n}$ be semisimple
$\Gamma$-modules.  Then $$V_{1}\tensor \ldots \tensor V_{n} \text{  is
semisimple if } p>\sum_{i=1}^{n}(\dim V_{i}-1).$$
\label{generalization:of:chevalleys:theorem}
\end{theorem}

The proof again  uses a reduction to algebraic
group theory.  As above we may assume that $\k$ is
algebraically closed, the representation $\Gamma \to \GL(V)$ is
faithful and  $\Gamma$ is a closed subgroup of $\GL(V)$ in the
Zariski topology, where $V=V_{1}\oplus \ldots \oplus V_{n}$.  But we
can no longer reduce to the case that $\Gamma$ is connected.  Indeed,
if $\Gamma$ is finite of order divisible by $p$, this assumption will
be no help at all.  So we need to do more. We need
$\Gamma$ to be \define{saturated}.

To define this notion (cf.  [N],[S1]), suppose that $x\in\GL_{n}(\k)$
has order $p$.  Write $x=1+\epsilon$ for some matrix $\epsilon$ and
note that $\epsilon^{p}=0$.  For any $t\in \k$  define
$x^{t}:=1+t\epsilon+\binom{t}{2} \epsilon^{2}+\ldots + \binom{t}{p-1}
\epsilon^{p-1}$.  Since $\epsilon^p = 0$ we have constructed a one
parameter subgroup $\{x^t \st t \in \k\}$ of $\GL_n(\k)$.  By
definition, a subgroup $\Gamma\subgroup\GL_n(\k)$ is said to be
\define{saturated} if it is Zariski closed and $x\in\Gamma$ with
$x^{p}=1$ implies that $x^{t}\in\Gamma$ for all $t\in\k$.  One can
define the \define{saturated closure} of a subgroup $\Gamma$ denoted
by $\Gamma\sat$.  It is the smallest saturated subgroup of
$\GL_{n}(\k)$ containing $\Gamma$.

Here are some examples:
\begin{itemize}
    \setlength{\itemsep}{0in}
    \item If $p>2$ every classical group in its natural
    representation is saturated.
    \item If $p>3$ the group $G_{2}(\k)$, embedded in $\GL_{7}(\k)$, is
    saturated.
    \item If $p=2$ the group $\PGL_{2}(\k)$, embedded in $\GL_{3}(\k)$ by its
    adjoint representation, is not saturated.
    \item If $p=11$ the Janko group $J_{1}$, embedded in $\GL_{7}(\k)$,
    has for saturated closure the group $G_{2}(\k)$.
    \end{itemize}

It can be checked that our problem is stable under replacing $\Gamma$
by $\Gamma\sat$.  So, we may assume that $\Gamma$ is saturated.  This
 implies that $\Gamma/\Gamma^{\circ}$ is finite of order prime to
$p$, so  we can reduce as before to the case where $\Gamma$ is a
connected reductive algebraic group.  Then we resort to the general
theory of representations of algebraic groups to complete the proof,
which is somewhat technical.  (cf. [S1])

One can also ask about various converse theorems (cf. [S2]).
For instance:
\begin{theorem_count}
    \MakeArabic
    \item Does $V_{1}\tensor V_{2}$ semisimple  imply  $V_{2}$
    semisimple?
    \item Does $\bigwedge^{2}V$ semisimple imply
    $V$   semisimple?
    \item Does $\Sym^{2} V$ semisimple imply   $V$   semisimple?
\end{theorem_count}
For question (1) the answer in
characteristic zero is yes unless $\dim V_{1}=0$.  In characteristic $p>0$,
the answer is yes unless $\dim
V_{1}=0$ in $\k$, i.e. unless $\dim V_{1}\equiv 0 \pmod{p}$.

For question (2) the answer in characteristic zero is yes
unless $\dim V=2$.  In characteristic $p>0$ the answer is yes unless
$\dim V\equiv 2 \pmod{p}$.

For question (3) the answer is yes
in characteristic zero, while in characteristic $p>0$ the answer is yes
unless $\dim V\equiv -2 \pmod{p}$.

\begin{remarks}
These questions make sense more generally in the setting of a
``tensor category'', cf. [D].  Such a category has
tensor products and duals, as well as a distinguished object
$\underline{1}$.  There is the notion of dimension of an object:
consider the composition of the natural maps
$$\underline{1}\to V\tensor V^{*}\to \underline{1}.$$
This  determines an element of $\k=\End(\underline{1})$, which
is called the dimension of $V$.  In particular it is
possible for the dimension to be $-2$ in $\k$. In this formalism,
there is a way of transforming symmetric powers into exterior powers,
by changing categories.  Deligne noticed that if one proves
in this setting one of the two statements:
\begin{center}
    $\bigwedge^{2}V$  semisimple $\Rightarrow$  $V$  semisimple
    if $\dim V\ne 2$ in $\k$\\
    $\Sym^{2}V$  semisimple $\Rightarrow$  $V$  semisimple if $\dim
    V\ne -2$ in $\k$
    \end{center}
then the other is true as well (cf. [S2],\S6.2).  (Here $\k$ is assumed to
be of
characteristic not equal to 2.)

W. Feit has provided various counterexamples showing that the results
are essentially the `best possible' for questions (1) and (2), (cf.
[S2], appendix).  The situation is different for question (3).  For
instance, with $p=7$ there is no known example in which $\Sym^{2}V$ is
semisimple but $V$ is not.
\end{remarks}

We turn now to giving a generalization of the notion of complete
reducibility (cf. [T2]).  Let $\k$ be algebraically closed, $G$ be a connected,
reductive algebraic $\k$-group and $\Gamma\subgroup G(\k)$.  I shall
say that 
$\Gamma$ is \define{$G$-completely reducible} ($G$-cr for short) if for every
parabolic subgroup $P$ of $G(\k)$ containing $\Gamma$ there exists a
Levi subgroup of $P$, also containing $\Gamma$.

The definition of $G$-cr may be reformulated within the context of Tits
buildings (cf. [T1]).  The Tits building of $G$
is the simplicial complex $X$, with simplices corresponding to the parabolic
subgroups of $G(\k)$ and inclusions being reversed.
The group $G(\k)$ acts simplicially on $X$. So if $\Gamma \subgroup G(\k)$,
we can consider the complex $X^\Gamma$ of all $\Gamma$-fixed points.
One can prove that there are precisely two possibilities:
\begin{theorem_count}
    \MakeArabic
\item  $X^\Gamma$ is contractible (homotopy type of a point);
\item  $X^\Gamma$ has the homotopy type of a bouquet of spheres.
\end{theorem_count}
One can show that (2) occurs precisely
when $\Gamma$ is $G$-cr.

The property of $\Gamma$ being $G$-cr relates nicely to the usual
property of a $\Gamma$-module being semisimple.  If we take $G$ to be
$\GL(V)$ for some vector space $V$, it is clear that $\Gamma$
is $G$-cr if and only if $V$ is a semisimple $\Gamma$-module.  More
generally, if $p\ne 2$ and $G$ is any symplectic group, orthogonal group,
or $G_{2}$ then $\Gamma$ is $G$-cr if and only if the natural
representation of $G(\k)$ is a semisimple $\Gamma$-module.
We would
like in a general setting, given $\Gamma \subgroup G(\k)$ and a linear
representation $V$ of $G(\k)$, to relate the property ``$\Gamma$ is
$G$-cr'' to the property that $V$ is a semisimple $\Gamma$-module (for
$p$ larger that some bound $n(V)$).  This will be discussed in the
later lectures.

Finally, we give an application of these ideas.  The Dynkin diagram of
$D_{4}$ has a symmetry of order 3 which gives rise to an automorphism
$\tau$ of $\Spin_8$.  Consequently, there are three irreducible
modules for $\Spin_{8}$ of dimension 8, say $V_{1}$, $V_{2}$, and
$V_{3}$.  Suppose that $\Gamma$ is
a subgroup of $\Spin_8$.  Is it true that:
$$V_1 \text{ is } \Gamma\text{-semisimple}\Rightarrow V_2 \text{ and }
V_3 \text{ are } \Gamma\text{-semisimple}?$$
The answer is yes if $p>2$ (and sometimes no if $p=2$): this follows
from the fact that $V_i$ is $\Gamma$-semisimple if and only if
$\Gamma$ is $\Spin_8$-cr.

\section*{Lecture 2}

Fix an algebraically closed field $\k$ and
let $G$ be a connected, reductive algebraic $\k$-group.
We are interested only in
the case where $p = \charistic \k>0$.  Recall that a subgroup
$\Gamma\subgroup G(\k)$ is called $G$-cr if for every parabolic
subgroup $P$ of $G(\k)$ containing $\Gamma$, there exists a Levi
subgroup of $P$ also containing $\Gamma$.  We wish to relate this to
the usual notion of complete reducibility. 

Let $T$ be a maximal torus of $G$, and $B$ be a Borel subgroup containing
$T$ with $U$ its unipotent radical.  This determines a root system and a
set of positive roots.  Let $X(T)=\Hom(T,\G_{m})$ be the character
group, and $Y(T)=\Hom(\G_{m},T)=\Hom(X(T),\Z)$ the cocharacter
group.  We have a natural pairing $\inprod{.,.}:X(T)\times Y(T)\to
\Z$ and for each $\alpha$ in the root system we have the coroot
$\alpha\check\in Y(T)$.

For each $\lambda \in X(T)$ define
$$
n_G(\lambda) = n(\lambda):= \sum_{\alpha>0}\inprod{\lambda,\alpha\check}.
$$
Note we can also write this as $\inprod{\lambda,\phi}$ where $\phi=
\sum_{ \alpha> 0} \alpha \check$, the {\em principal homomorphism} of
$\G_m$ into $T$.  If $V$ is any finite dimensional $G$-module, let us put;
$$
n_{G}(V) = n(V) :=\sup n(\lambda)
$$
where the supremum is taken over all the weights $\lambda$ of $T$ in  $V$.

As an example, consider $G = \GL_{m}$, with $V$ the natural
$m$-dimensional representation.
Then $n(V)=m-1=\dim V -1$ and $n(\bigwedge^{i}V) =
i(\dim V -i)$.
In general, if $V_1$ and $V_2$ are any $G$-modules,
$n(V_{1}\tensor V_{2}) = n(V_{1}) +
n(V_{2})$.

Note that if $V$ is a {\em nondegenerate} linear representation of
$G$, i.e.  the connected kernel of the representation is a torus,
then $n(V)\ge h-1$, where $h$ is the Coxeter number of $G$.  Indeed,
let
$\lambda$  be a highest weight of $V$.  So
$n(V)=n(\lambda)=\inprod{\lambda,\sum_{\alpha>0}\alpha\check} \ge
\inprod{\lambda+\rho,\beta\check}-1$, where $\rho$ is half the sum
of positive roots and $\beta\check$ is the highest coroot.  Since
$\lambda$ is nonzero and dominant we have
$\inprod{\lambda,\beta\check}\ge 1$  and
$\inprod{\rho,\beta\check}=h-1$.

Our goal is to prove the following result:

\begin{main_theorem}  Let $V$ be $G$-module with $p>n(V)$.  Let
$\Gamma$ be a subgroup of $G(\k)$.  Then
$$\Gamma \text{ is } G\text{-cr} \Rightarrow V \text{ is }
\Gamma\text{-semisimple}.$$
Moreover, the converse is true if $V$ is nondegenerate, i.e. the connected
kernel of the representation is a torus.
\end{main_theorem}

Some of the results mentioned in Lecture 1 are immediate
consequences.  For example, let $\{V_{1}, \dots, V_m\}$ be a
collection of semisimple $\Gamma$-modules and $p>\sum_i (\dim
V_{i}-1)$.  Then the theorem, applied to $G=\prod\GL(V_{i})$ and
$V=\bigotimes V_{i}$, tells us that $V_1 \otimes \dots \otimes
V_m$ is also semisimple.  Alternatively, suppose that $V$ is
a semisimple $\Gamma$-module with $p>i(\dim V-i)$.  Then the
theorem shows that $\bigwedge^{i}V$ is $\Gamma$
semisimple.  (This was stated as an open question at the end of [S2];
and the special case where $V$ is irreducible had been proved by
McNinch.)\bigskip

The proof of the main theorem uses the notion of
\define{saturation} with respect to the group $G$. In order to define
it, we need to introduce the ``exponential'' $x^{t}$, for $x$
unipotent in $G$ and $t\in\k$.
This is possible (for $p$ not too small) thanks to:

\begin{theorem}
Assume $p \geq h$ (resp. $p > h$ if $G$ is not simply
connected).
There exists a unique isomorphism of
varieties $\log:G^{u}\to\g_{\nilp}$ with the following
properties:
\begin{theorem_count}
    \item $\log(\sigma u)=\sigma\log u$ for all $\sigma \in \Aut G$;
    \item the restriction of $\log$ to $U(\k)$ is an isomorphism of
    algebraic groups $U(\k) \to \Lie U$, whose tangent
    map is the identity;
    \item $\log(x_{\alpha}(\theta))=\theta X_{\alpha}$, for every
    root $\alpha$ and every $\theta\in\k$.
    \end{theorem_count}
    \label{existence:of:log}
\end{theorem}

Here, $\g_{\nilp}$ is the nilpotent variety of $\Lie G$,
$x_{\alpha}:\G_a \to U_\alpha$ denotes some fixed parameterization of
the root group $U_\alpha$ of $U$, and $X_{\alpha} =
\frac{d}{d\theta}(x_\alpha(\theta))|_{\theta = 0}$ is the
corresponding basis element of $\Lie U_\alpha$.  We are viewing $\Lie
U$ as an algebraic group over $\k$ via the Campbell-Hausdorff formula:
$XY:=X+Y+\frac{1}{2}[X,Y]+\frac{1}{12}[X,[X,Y]]+\ldots$ (cf.  [B],
Chap II, \S6) which makes sense in characteristic $p$ because of the
assumption $p\ge h$ and the fact that the nilpotency class of $\Lie U$
is at most $h$.

For the proof, uniqueness is obvious since the $U_{\alpha}$ generate
$U$ and $G^{u}$ is the union of conjugates of $U$.  (Moreover, one can
show that (iii) is a consequence of (i) and (ii).)  However, the
existence part is less easy.  One possible method is to define first
$\log$ on $U$ and then extend it to $G^{u}$.  This approach uses the
fact that $\g_{\nilp}$ is a normal variety (cf. [D], [BR]) when $p$
is good, that $G^{u}$ is a normal variety (cf. [St]) and draws on work
by Springer (cf. [Sp2]).

Given the theorem,
let $\exp:\g_{\nilp} \to G^u$ denote the inverse to $\log$.
For  $x\in G^{u}(\k)$ and $t\in \k$ we
define $x^{t}$ as $\exp(t\log x)$.
Note that the exponential map $x,t\mapsto x^{t}$ may be viewed as a
morphism $F:G^{u}\times \A^{1}\to G^{u}$.  Moreover this map is the
``reduction $\operatorname{mod} p$''
of the corresponding  well-known map in
characteristic zero, and this gives a convenient way to compute it.


\section*{Lecture 3}

Continue with the notation of the previous lecture.  Recall that we
have just defined the map $x \mapsto x^t$ for any unipotent element $x
\in G(\k)$ and any $t \in k$.  We can now at last define the saturation
process (assuming $p\ge h$).  A subgroup $\Gamma$ of $G(\k)$ is \emph{saturated} if
\begin{theorem_count}
    \MakeArabic
    \item  $\Gamma$ is Zariski closed;
    \item whenever $x\in \Gamma\cap G^{u}$, we have
$x^{t}\in \Gamma$ for all $t \in \k$.
\end{theorem_count}
We wish in the remainder of this lecture to describe some basic
properties of saturated subgroups and $G$-cr subgroups.  We will apply
these properties in Lecture 4 to prove the Main Theorem.

We begin by mentioning some elementary examples: every parabolic
subgroup is saturated; the centralizer of any subgroup of $G(\k)$ is
saturated; Levi subgroups are saturated, since they may be realized as
the centralizer of a torus.  We also note that in the case of
saturated subgroups lying in $U$, there are various alternative
characterizations giving further `unipotent' examples:

\begin{property}
Let $V$ be a closed subgroup of $U(\k)$.  The following are equivalent:
\begin{theorem_count}
    \item  $V$ is saturated;
    \item  $V=\exp(\mathfrak{v})$ for $\mathfrak{v}$ a Lie subalgebra
    of $\Lie U$;
    \item  $\log V$ is a vector subspace of $\Lie U$.
    \end{theorem_count}
    \label{characterize:saturated}
\end{property}

Another basic property is as follows:

\begin{property}
Let $H$ be a semisimple subgroup of $G$ with $H(\k)$ saturated.  If $x$
is any unipotent element of $H(\k)$, then the element $x^t$ (defined
relative to $H$) coincides with $x^t$ (defined relative to $G$).
\label{subgroups:preserve:exp}
\end{property}

Even to state Property \ref{subgroups:preserve:exp} correctly, we need
first to know that the
Coxeter number $h_H$ of $H$ does not exceed the Coxeter number $h_G$
of $G$.  In fact, a stronger result holds:

\begin{theorem}
Let $p$ be any prime, and $H$ be a semisimple subgroup of a semisimple
group $G$.  Let $d_{i, H}$ and $d_{j, G}$ be the invariant degrees of
the Weyl groups of $H$ and $G$ respectively.  Then, the polynomial
$\prod(1-T^{d_{i, H}})$ divides $\prod(1-T^{d_{j, G}})$.\newline
{\rm (For the properties of the invariant degrees, see [B], Chap V,
\S5.)}
\label{theorem:poincare:series}
\end{theorem}

As a corollary we see that each $d_{i, H}$ divides some $d_{j, G}$.
For, choosing $T$ to be a primitive $d_{j, H}\th$ root of unity, the
theorem implies that $\prod(1-T^{d_{j, G}})$ vanishes, so $(1 -
T^{d_{j, G}})$ vanishes for some $j$.  Since the Coxeter number $h_H$
is the largest degree $d_{i, H}$, and similarly for $G$, we deduce in
particular that $h_{H}\le h_{G}$, as required for the statement of
Property \ref{subgroups:preserve:exp} to make sense.



We sketch the proof of Property \ref{subgroups:preserve:exp}.
Assume that $H$ is a semisimple
subgroup of $G$ with $H(\k)$ saturated.  We may assume that there is a
maximal unipotent subgroup $U_H$ of $H$ with $U_H \subgroup U$.  Note
that $U_H(\k)$ is also a saturated subgroup of $G$.  We need to show
that $\log_{G}(x)=\log_{H}(x)$ for any unipotent $x\in H(\k)$.
Conjugating, it suffices to prove this for $x \in U_H(\k)$.  We have an
isomorphism $\log:U(\k)\to\Lie U$.  Viewing $\Lie U_{H}$ as a
subgroup of $\Lie U$, we conclude that the restriction of $\log_G$
gives a isomorphism $U_{H}\cong \Lie U_{H}$ which is compatible with
conjugation and whose tangent map is the identity.  By the uniqueness
in the definition of $\log_H$ we conclude that the restriction of
$\log_G$ is equal to $\log_H$, as required.

\begin{property} If $H\subgroup G$ is saturated then the index
$\SGindex{H}{H^{\circ}}$ is prime to $p$.
\label{index:prime}
\end{property}
To prove Property 3, suppose $p$ divides $\SGindex{H}{H^{0}}$ and take some
element $x$ of the finite group $H/H^{0}$ of order $p$. One proves,
from general principles, that there exists  $\tilde{x}\in H^{u}(\k)$ which
maps
onto $x$ in the quotient.  By saturation, $\{\tilde x^t \st t \in k\}$
is a subgroup of $H(\k)$, hence of $H^0(\k)$ since it is connected.  So
$\tilde x\in H^0(\k)$, a contradiction.

We turn to discussing some basic properties of $G$-cr subgroups, as
defined in Lectures 1 and 2.  Recall that given a completely reducible
$H$-module for an algebraic group $H$, the unipotent radical of $H$
acts trivially.  The next property that we will need is similar, but
stated intrinsically within the groups.

\begin{property}
    If $\Gamma$ is $G$-cr and $V$ is a normal unipotent subgroup of
    $\Gamma$ then $V=1$.  In particular, if in addition $\Gamma$ is
    Zariski closed, then $\Gamma^{0}$ is reductive.
    \label{Gcr:implies:reductive}
\end{property}

The proof of this depends on the construction of Borel and Tits (cf.
[BT]) which associates to the subgroup $V$  a
parabolic subgroup $P$ of $G$ with $V\subgroup \urad(P)$.  Now
$\Gamma$ normalizes $V$, and since $P$ is defined in a canonical
fashion, $\Gamma$ normalizes $P$.  Therefore
$\Gamma\subgroup P$.  Now we use the fact that $\Gamma$ is $G$-cr to
deduce that $\Gamma\cap \urad(P)=1$, whence $V=1$.

\begin{property}
    Let $\Gamma_{0}\subgroup \Gamma$ be a normal subgroup of $\Gamma$
    of index prime to $p$.  Then, $\Gamma_{0}$ is $G$-cr $\Rightarrow
    \Gamma$ is $G$-cr.
    \label{conn:comp:Gcr:implies:Gcr}
\end{property}

(Before sketching the proof of this, we mention an open problem: if
$\Gamma_{0}\subgroup \Gamma$ is normal, is it true that $\Gamma$ is
$G$-cr $\Rightarrow \Gamma_{0}$ is $G$-cr?)

Now for the proof, let $P$ be a parabolic subgroup of $G$
containing $\Gamma$, and let $L$ be a Levi subgroup of $P$ which
contains $\Gamma_{0}$.  Write $P=\urad P\semidirect L$.  Let
$\Gamma_{L}$ be the image of $\Gamma$ under the projection $P \to L$.
The kernel of this projection is $\Gamma\cap\urad P=1$ so we have an
isomorphism $\Gamma\to \Gamma_{L}$.  Then $\Gamma$ is obtained from
$\Gamma_{L}$ by a 1-cocycle $a:\Gamma\to \urad P$, with $a$ equal to a
coboundary on restriction to $\Gamma_{0}$.  This implies that $a$ is
induced by a 1-cocycle $a'$ on $\Gamma/\Gamma_{0}$ with values in
$V=\urad P\cap Z(\Gamma_{0})=(\urad P)^{\Gamma_{0}}$.  Now, $V$ has a
composition series made up of $\k$-vector spaces, and since
$\size{\Gamma/\Gamma_{0}}$ is prime to $p$, the cocycle induced by $a'$
on each such composition factor is a coboundary.  This implies that
$a'$, whence $a$, is a coboundary, so that we can conjugate $\Gamma$
to a subgroup of $L$, as required.

\section*{Lecture 4}

Now we proceed to prove the Main Theorem. We begin with:

\begin{theorem}
    Suppose $p\ge h$.  Let $V$ be a $G$-module with associated
    representation $\rho_V:G \to \GL(V)$.  For every unipotent
    element $u$ of $G$, let $d_{u}(V)$ be the degree of the polynomial
    map $t\mapsto \rho_{V}(u^{t})\in \End(V)$.  Then $d_{u}(V) \le
    n(V)$, and there is equality if $u$ is regular.
    \label{compare:degrees:invariant}
\end{theorem}

The proof is in several steps.

(1)  \emph{The case $G=\SL_{2}$}.  In this case we may assume $u=\twobytwo
1101$, $u^{t}=\twobytwo 1t01$, and we have to prove $d_{u}(V) = n(V)$.

(1.1) \emph{One has $d_{u}(V) \le n(V)$}.  Write $\rho_{V}(u^{t})$ as
$1+\sum_{i\ge 1}a_{i}t^{i}$, $a_{i}\in \End(V)$.  If $s_{\lambda} =
\twobytwo {\lambda}00{\lambda^{-1}}$ with $\lambda\in \k^{*}$, we
have $s_{\lambda}u^{t}s_{\lambda}^{-1} = u^{\lambda^{2}t}$, hence
$$
\rho^{V}(s_{\lambda}).\sum a_{i}t^{i}.\rho_{V}(s_{\lambda}^{-1}) =
\sum a_{i}\lambda^{2i}t^{i},
$$
which implies $\rho_{V}(s_{\lambda}) a_{i}
\rho_{V}(s_{\lambda})^{-1} = \lambda^{2i}a_{i}$ for every $i$.  Hence
$a_{i}$ has weight $2i$ in $\End(V)= V\tensor V^{*}$.  By definition
of the invariant $n$ this shows that $a_{i}\ne 0 \Rightarrow 2i\le
n(V\tensor V^{*}) = n(V)+n(V^{*})=2n(V)$, i.e. $i\le n(V)$.  Hence
$d_{u}(V) \le n(V)$.

(1.2) \emph{One has $d_{u}(V)\ge n(V)$}.  If $V$ has Jordan-H\"older
quotients $V_{\alpha}$, it is clear that $n(V) =\sup n(V_{\alpha})$,
$d_{u}(V)\ge \sup d_{u}(V_{\alpha})$.  Hence we may assume that $V$
is simple.  In that case, the equality $n(V)=d_{u}(V)$ is obvious
from the explicit description of $V$ \`a la Steinberg.

(2)  \emph{The case $G$ arbitrary, $u$ regular}.  Choose a principal
homomorphism $\SL_{2}\to G$,  (cf. [Te] - see also
[S3]).  It is known that a nontrivial unipotent of
$\SL_{2}$ gives a regular unipotent of $G$.  On the other hand, one
has $n_{G}(V) = n_{\SL_{2}}(V)$, almost by definition.  Hence the
result follows from (1).

(3) \emph{General case}.  For $u$ unipotent of $G$, write
$\rho_{V}(u^{t})$ as above:
$$
\rho_{V}(u^{t})=1+\sum a_{i}(u)t^{i}\in \End(V).
$$ The $a_{i}$ are regular functions of $u$ (viewed as a point of the
unipotent variety $G^{u}$).  If $i<n(V)$ then  $a_{i}(u)$ is 0 when $u$ is
regular, by (2).  Since the regular unipotents are dense in $G^{u}$,
this implies $a_{i}(u)=0$ for every $u$.

\begin{cor}
    If $H$ is a reductive and saturated subgroup of $G$, one has
    $n_{H}(V) \le n_{G}(V)$.
    \label{invariant:decreases:subgroups}
\end{cor}

Choose a regular unipotent element $u\in H$.  One gets $n_{H}(V) =
d_{u}(V) \le n_{G}(V)$ by Theorem \ref{compare:degrees:invariant}, applied to both $H$ and $G$.

\begin{cor}
    The following are equivalent:
    \begin{theorem_count}
	\item $p>n(V)$;
	\item $\rho_{V}(u^{t}) = \rho_{V}(u)^{t}$ for every unipotent $u$ of
	$G$, and every $t\in \k$.
    \end{theorem_count}
    \label{invariance:of:exp:through:rho}
\end{cor}

Indeed (ii) holds if and only if the degree of $t\mapsto
\rho_{V}(u^{t})$ is $<p$, i.e. if and only if $d_{u}(V)<p$.
Since $n(V)=\sup_{u}d_{u}(V)$, this shows the equivalence of (i)
and (ii).  (The same proof shows that (i) and (ii) are equivalent to:

\emph{{\rm(ii$'$)} $\rho_{V}(u^{t})=\rho_{V}(u)^{t}$ for every regular $u$,
and every $t\in \k$.})

\begin{theorem}
    Let $G$ be reductive connected, and let $V$ be a $G$-module.
    Assume $p>n(V)$.  Let $\Gamma$ be a subgroup of $G(\k)$, which is
    $G$-cr.  Then $V$ is $\Gamma$-semisimple.
    \label{Gcr:implies:ss}
\end{theorem}

The proof is in several steps.

(1) We may assume that $\rho_{V}:G\to \GL(V)$ has trivial kernel.

(2) \emph{We have $p\ge h$}.  This follows from $p>n(V)\ge h-1$ (cf.
Lecture 2).

(3) \emph{The $G$-module $V$ is semisimple}.  Write $G$ as
$T.S_{1}\ldots S_{m}$, where $T$ is the maximal central torus, and
$S_{1}\ldots S_{m}$ is the decomposition of $(G,G)$ into quasi-simple
factors.  To prove (3), it is enough to show that $V$ is
$S_{i}$-semisimple for every $i$ (this is an easy lemma, cf. [J2]
and comments below); since $n_{s_{i}}(V) \le n_{G}(V)$ we are reduced
to the case where $G$ is quasi-simple.  With the usual notation we
have, for every weight $\lambda$ of $V$, $\lambda\ne 0$,
$$
\inprod{\lambda+\rho,\beta\check}\le
1+\sum_{\alpha>0}\inprod{\lambda,\alpha\check} \le p,
$$
where the inequality on the left is in [S1], p.519.  This shows that
the simple modules $L(\lambda_{i})$ in a Jordan-H\"older series of $V$
are of two types: $\lambda_{i}=0$, or
$\inprod{\lambda_{i}+\rho,\beta\check}\le p$.  But it is known (cf.
[J1]) that this implies $L(\lambda_{i})$
$\Ext_{G}^{1}(L(\lambda_{i}),L(\lambda_{i}))=0$ for every pair
$\lambda_{i},\lambda_{j}$ (e.g. because these $L(\lambda_{i})$ are
Weyl modules).  Hence $V$ is semisimple.

We pause to discuss a variant of this proof.  If $\lambda$ is a
dominant weight with
$\sum_{\alpha>0}\inprod{\lambda,\alpha\check}<p$, then $L(\lambda) =
V(\lambda)$, where $V(\lambda)$ is the Weyl module.  The proof is by
reduction to $G$ quasi-simple, and one distinguishes between two
cases: $\lambda=0$, where it is obvious, and $\lambda\ne 0$,  where we
have
$\inprod{\lambda+\rho,\beta\check} \le p$. Moreover, if $\lambda$, $\mu$
have the property $L(\lambda)=V(\lambda)$, $L(\mu) = V(\mu)$ then
$\Ext^{1}_{G}(L(\lambda),L(\mu))=0$.  See [J2] for more details.

(4) \emph{The $\Gamma\sat$-module $V$ is semisimple}.  (Note that
we can define $\Gamma\sat$  since $p\ge h$ by (2).)  Let $H$ be the
connected component of $\Gamma\sat$.  Since $\Gamma\sat$ is $G$-cr
(because $\Gamma$ is), $H$ is a reductive group.  By Corollary
\ref{invariant:decreases:subgroups} to
Theorem \ref{compare:degrees:invariant}, we have $n_{H}(V) \le n_{G}(V)$
hence $n_{H}(V)<p$ and
part (3) above (applied to $H$) shows that $V$ is $H$-semisimple.
Since $\SGindex{\Gamma\sat}{H}$ is prime to $p$, this implies that $V$ is
$\Gamma\sat$-semisimple (cf. [S1], p.523).

(5)  \emph{If a subspace $W$ of $V$ is $\Gamma$-stable, it is
$\Gamma\sat$-stable}.  Let $H_{W}$ be the stabilizer of $W$ in $G$.
If $u\in H_{W}$ is unipotent, one has $\rho_{V}(u^{t}) = \rho_{V}
(u)^{t}$ by Corollary \ref{invariance:of:exp:through:rho}
to Theorem \ref{compare:degrees:invariant} above.  Since $\rho_{V}(u)W=W$
the same is true for $\rho_{V}(u)^{t}$ for every $t$.  This shows that
$H_{W}$ is saturated.  Since it contains $\Gamma$, it also contains
$\Gamma\sat$.

(6) \emph{End of proof}.  By (5), the subspaces of $V$ which are
$\Gamma$-stable are the same as those which are $\Gamma\sat$-stable.
Since, by (4), $V$ is $\Gamma\sat$-semisimple, it is
$\Gamma$-semisimple.

Note that this is the ``Main Theorem'' announced at the beginning of
these lectures.  It implies for instance the following (where $\k$ is
arbitrary of characteristic $p$):

\emph{If $V_{\alpha}$ are semisimple $\Gamma$ modules, and
$i_{\alpha}\ge0$  integers with
$$\sum i_{\alpha}(\dim
V_{\alpha}-i_{\alpha})<p,$$
 then
$\bigotimes_{\alpha}\bigwedge^{i_{\alpha}}V_{\alpha}$ is semisimple.}

(Sketch of proof.  Apply Theorem \ref{Gcr:implies:ss} to
$\prod_{\alpha}\GL(V_{\alpha})$
and $V=\bigotimes_{\alpha}\bigwedge^{i_{\alpha}}V_{\alpha}$, and
deduce the  statement when $\k$ is algebraically closed.  Next show that
one can assume $i_{\alpha} \le (\dim V_{\alpha})/2$, and $\dim
V_{\alpha} <p$ for all $\alpha$; deduce that $V_{\alpha}$ is
absolutely semisimple (i.e. remains semisimple after extension of
scalars from $\k$ to $\overline{\k}$); hence
$\bigotimes_{\alpha}\bigwedge^{i_{\alpha}}V_{\alpha}$ is absolutely
semisimple.)

\begin{theorem}[Eugene] {\rm (cf. [J2], [Mc], [LS])}
    Let $H\subgroup G$ be connected reductive, and $p\ge h_{G}$.  Then
    $H$ is $G$-cr.
    \label{conn:red:implies:Gcr}
\end{theorem}

The proof starts by reducing to the case $G$ is quasi-simple.  Then
there are separate proofs for type $A_{n}$ (Jantzen), $B_{n}$,
$C_{n}$, $D_{n}$ (Jantzen-McNinch), and exceptional type
(Liebeck-Seitz).  There is a little extra work involved in the
$B_{n}$, $C_{n}$, $D_{n}$ cases when $H$ is of type $A_{1}$.  (Note
that in special cases $p\ge h_{G}$ can be improved.)

\begin{theorem}
    Let $\Gamma\subgroup G$.  Assume $p\ge h_{G}$.  The following are
    equivalent:
    \begin{theorem_count}
	\item $\Gamma$ is $G$-cr;
	\item the connected component of $\Gamma\sat$ is reductive.
	\end{theorem_count}
\label{Gcr:equiv:conn:comp:red}
\end{theorem}

The direction (i)$\Rightarrow$(ii) is clear since $\Gamma$ is $G$-cr
$\Rightarrow$ $\Gamma\sat$ is $G$-cr, and hence $(\Gamma\sat)^{0}$ is
reductive.

For (ii)$\Rightarrow$(i) apply Theorem \ref{conn:red:implies:Gcr}
to $H=(\Gamma\sat)^{0}$.  One
sees that $H$ is $G$-cr, hence also $\Gamma\sat$, hence also
$\Gamma$.

\begin{theorem}
    Let $V$ be a nondegenerate $G$-module.  Assume $n(V)<p$.  If
    $\Gamma\subgroup G$, the following are equivalent:
    \begin{theorem_count}
	\item  $\Gamma$ is $G$-cr;
	\item $V$ is $\Gamma$-semisimple.
	\end{theorem_count}
\label{Gcr:equiv:ss}
\end{theorem}

The direction (i)$\Rightarrow$(ii) is Theorem \ref{Gcr:implies:ss}.
Conversely, if $V$
is $\Gamma$-semisimple, it is also $\Gamma\sat$-semisimple (cf.
argument of Theorem 6), hence $(\Gamma\sat)^{0}$-semisimple and by
Theorem 8 this shows that $\Gamma$ is $G$-cr.

Note: The implication (ii)$\Rightarrow$(i) proved above under the
condition $p>n(V)$ is far from best possible.  Example: take
$G=\GL(W)$, and $V=\bigwedge^{2}W$, which is nondegenerate if $\dim
W\ne 2$.  One has $n(V)=2(\dim W-2)$ and one sees that:
$$
{\textstyle\bigwedge^{2}}W \text{ is } \Gamma\text{-semisimple}\ \Rightarrow\ W
\text{ is } \Gamma\text{-semisimple}
$$
if $p>2(\dim W-2)$.  However, an elementary argument [S2],
shows that this remains true as long as $p$ does not
divide $\dim W-2$.

Example of Theorem 8:  If one takes for $V$ the adjoint
representation $\Lie G$, which is nondegenerate, one has $n(\Lie
G)=2h-2$ and Theorem 8 gives:
$$
\Gamma \text{ is } G\text{-cr} \iff \Lie G \text{ is }
\Gamma\text{-semisimple}
$$
provided $p>2h-2$.  (In fact, for $G=\GL_{n}$, no condition on $p$ is
needed for $\Leftarrow$, cf. [S2], Theorem
3.3.)
\newpage

\thispagestyle{empty}
\vspace*{20mm}

\begin{center}
{\bf\Huge A few references}
\end{center}
\vspace{5mm}

\begin{itemize}
\item[{[BR]}] P. Bardsley and R. W. Richardson,
\'Etale slices for algebraic transformation groups in characteristic
$p$, {\em Proc. London Math. Soc.} {\bf 51} (1985), 295--317.
\item[{[Bo1]}]  A. Borel, Sous-groupes commutatifs et torsion des
groupes de Lie compacts connexes, {\em T\^ohoku Math. J.}, (2) \textbf{13}
(1961), 216-240
\item[{[BS]}]
A. Borel and J.-P. Serre, Sur certains sous-groupes des groupes de
Lie compacts, {\em Comm. Math. Helv.} {\bf 27} (1953), 128--139.
\item[{[BT]}] A. Borel and J. Tits, \'El\'ements unipotents et
sous-groupes paraboliques des groupes r\'eductifs, \emph{Invent.
math.} \textbf{12} (1971), 95-104.
\item[{[B]}] N. Bourbaki, {\em Groupes et alg\`ebres de Lie}, Chap I-VI,
Hermann, Paris (1968).
\item[{[C]}]
C. Chevalley, {\em Th\'eorie des groupes de Lie}, vol. III,
Hermann, Paris, (1954).
\item[{[D]}] P. Deligne, Cat\'egories Tannakiennes, \emph{The
Grothendieck Festschrift}, Vol II, Birkh\"auser, Boston (1990), 111-195.
\item[{[De]}] M. Demazure, Invariants sym\'etriques entiers des
groupes de Weyl et torsion, \emph{Invent. math.} \textbf{21} (1973),
287-301.
\item[{[Dy]}] E. Dynkin, Semisimple subalgebras of semisimple Lie
algebras, \emph{A.M.S. Translations series 2} \textbf{6} (1957),
111-245.
\item[{[G]}]
R. Griess,
Elementary abelian $p$-subgroups of algebraic groups,
{\em Geom. Dedicata} {\bf 39} (1991), 253--305.
\item[{[GR]}]
R. Griess and A. Ryba, Finite simple groups which projectively embed
in an exceptional Lie group are classified!, preprint, April 14, (1998).
\item[{[J1]}] J. Jantzen, {\em Representations of Algebraic Groups},
Academic Press, Orlando, Pure and Applied Mathematics \textbf{131} (1987).
\item[{[J2]}] J. Jantzen, Low dimensional representations of reductive
groups are semisimple, \emph{Algebraic Groups and Lie groups: a
volume in honor of R.W. Richardson.} Cambridge (1997).
\item[{[K]}] B. Kostant, The principal three-dimensional subgroup and
the Betti numbers of  a complex simple Lie group, \emph{Amer. J. Math.}
\textbf{81} (1959), 973-1032.
\item[{[LS]}] M. Liebeck G. Seitz, Reductive subgroups of exceptional
algebraic groups, \emph{Memoirs of the A.M.S.}
\textbf{580} (1996).
\item[{[Mc]}] G. McNinch, Dimensional criteria for semisimplicity of
representations, \emph{Proc. London
Math. Soc.}(3) \textbf{76} (1998), 95--149.
\item[{[Min]}] H. Minkowski, {\em Gesammelte Abhandlungen}, New York,
Chelsea Pub. Co. (1967).
\item[{[N]}]
M. V. Nori, On subgroups of $\GL_n(\F_q)$, {\em Invent. math.}
 {\bf 88} (1987), 257--275.
\item[{[S1]}]
J.-P. Serre, Sur la semi-simplicit\'e des produits tensoriels de
repr\'esenta\-tions de groupes, {\em Invent. math.} {\bf 116} (1994),
513--530.
\item[{[S2]}]
J.-P. Serre, Semisimplicity and tensor products of group representations:
converse theorems, {\em J. Algebra} {\bf 194} (1997), 496--520 (with an
appendix
by W. Feit).
\item[{[S3]}]
J.-P. Serre, Exemples de plongements des groupes $\PSL_2(\F_p)$
dans des groupes de Lie simples,
{\em Invent. math.} {\bf 124} (1996), 525--562.
\item[{[S4]}] J.-P. Serre, {\em Galois Cohomology}, Springer-Verlag (1997).
\item[{[Sp1]}] T. Springer, Some arithmetical results on semi-simple
Lie algbras, \emph{Publ. Math I.H.E.S.} \textbf{30}, (1966),115-141.
\item[{[Sp2]}]
T. Springer, Regular elements of finite reflection groups,
{\em Invent. math.} {\bf 25} (1974), 159--198.
\item[{[SpSt]}] T. Springer, R. Steinberg,  Seminar
on algebraic groups and related finite groups: conjugacy classes,
algebraic groups, \emph{Lect. Notes in Math} \textbf{131},
Springer-Verlag (1970).
\item[{[St]}]
R. Steinberg, Regular elements of semisimple algebraic groups,
{\em Publ. Math. IHES} {\bf 25} (1965), 49--80 (= Coll. Papers no. 20).
\item[{[Te]}] D. Testerman, The construction of the maximal $A_{1}$'s in
the exceptional algebraic groups. \emph{Proc. A.M.S} \textbf{161},
635-644 (1993).
\item[{[T1]}]
J. Tits, Buildings of spherical type and finite $BN$-pairs, \emph{Lect.
Notes in Math.} \textbf{386},
Springer-Verlag (1974).
\item[{[T2]}] J. Tits, R\'esum\'e de cours au Coll\`ege de France,
(1997-1998), 93-98.
\end{itemize}
\end{document}